\documentclass[letterpaper,11pt]{article}
\usepackage{bm}
\usepackage{xcolor}
\usepackage{graphicx}
\usepackage{anyfontsize}
\usepackage[round,authoryear]{natbib}
\usepackage{natbib}
\usepackage{url}
\usepackage{multirow}
\usepackage{amstext}
\usepackage{amssymb}
\usepackage{amsmath}
\usepackage{lmodern}
\usepackage{textcomp}
\usepackage{afterpage,lscape}
\usepackage{tikz}
\usepackage{smartdiagram}
\usepackage[utf8]{inputenc}
\usepackage[T1]{fontenc}
\usepackage{textcomp}
\usepackage{caption}
\usepackage{subcaption}
\usepackage{tikz}
\usepackage{multirow}
\usepackage{here}

\usepackage{indentfirst}
\allowdisplaybreaks
\usepackage{longtable}
\usepackage{enumitem}
\usepackage{booktabs}
\newtheorem{property}{Property}

\title{Workload Equity in Multi-Period Vehicle Routing Problems}

\author{Najmeh Nekooghadirli, Michel Gendreau, Jean-Yves Potvin, Thibaut Vidal}

\begin{document}

\linespread{1.25}\selectfont

\vspace*{0.2cm}

\begin{center}
\begin{Huge}
Workload Equity in Multi-Period\vspace*{0.3cm}\linebreak Vehicle Routing Problems
\end{Huge}

\vspace*{0.5cm}

\textbf{Najmeh Nekooghadirli$^{1}$, Michel Gendreau$^{1}$}\\ 
\textbf{Jean-Yves Potvin$^{2}$, Thibaut Vidal$^{1,3,4}$}\\ \vspace*{0.3cm}

{\small
$^{1}$ CIRRELT \& Department of Mathematical and Industrial Engineering, Polytechnique Montreal, Canada\\

$^{2}$ CIRRELT \& Department of Computer Science and Operation Research, University of Montreal, Canada\\ 

$^{3}$ Department of Computer Science, Pontifical Catholic University of Rio de Janeiro, Brazil\\

$^{4}$ SCALE-AI Chair in Data-Driven Supply Chains\\
}

\end{center}

\vfill

\noindent{\bf{Abstract.}}
An equitable distribution of workload is essential when deploying vehicle routing solutions in practice. For this reason, previous studies have formulated vehicle routing problems with workload-balance objectives or constraints, leading to trade-off solutions between routing costs and workload equity. These methods consider a single planning period; however, equity is often sought over several days in practice. In this work, we show that workload equity over multiple periods can be achieved without impact on transportation costs when the planning horizon is sufficiently large. To achieve this, we design a two-phase method to solve multi-period vehicle routing problems with workload balance. Firstly, our approach produces solutions with minimal distance for each period. Next, the resulting routes are allocated to drivers to obtain equitable workloads over the planning horizon. We conduct extensive numerical experiments to measure the performance of the proposed approach and the level of workload equity achieved for different planning-horizon lengths. For horizons of five days or more, we observe that near-optimal workload equity and optimal routing costs are jointly achievable.\\ 

\noindent{{\bf Keywords:} Vehicle routing problem, Multiple periods, Workload equity, Optimization}\\

\vfill
\noindent
$^*$ Correspondence to: \url{thibaut.vidal@polymtl.ca}

\pagebreak

\section{Introduction}

Competition between companies in the transportation and distribution sectors has created a need to improve many business practices. Although minimizing distribution costs is essential, drivers' and clients' satisfaction should also not be neglected. Equitable workload assignments, especially, are fundamental to maintaining employee satisfaction at the workplace and can significantly impact the quality of services provided to clients. In recent works on vehicle routing problems with equity considerations, workload balancing is often treated as a second objective---along with cost minimization---in a bi-objective problem. In this context, \cite{articleMatl2019} reported that the marginal cost of balancing routes is reasonable in most cases: nearly 40\% of Pareto-optimal solutions have an additional cost that does not exceed 10\% of the minimum-cost solution. Solutions that account for equity also appear more robust against unexpected events, such as a sudden increase in demand, because bottlenecks are reduced \citep{articleMatl2019, MOURGAYA20071028}. Yet, such extra costs remain significant given that the transportation sector operates with very tight profit margins.

In most previous studies, we noted that the definition of workload equity was unnecessarily restrictive, as it focused on the balance between different drivers' workloads for each day separately. In practical situations, however, workload differences among drivers can be acceptable on each given day as long as the total workload over a longer time horizon (e.g., a week) remains equitable. To capture this notion, we formally define and study a multi-period VRP with workload balance (MVRPB). In this problem, customer requests are known on a longer planning horizon, and the goal is to find routes that optimize distance and workload balance over the entire planning horizon. To solve the MVRPB, we propose a two-phase solution approach. In the first phase, a solution with optimal distance is found for each period by optimally solving the corresponding VRPs. The resulting routes are then combined in a second phase to obtain equitable workloads among drivers based on the total distance traveled over the periods. We use this methodology to evaluate to which extent the length of the planning horizon can impact equitable solutions. Therefore, this work makes the following contributions:
\begin{itemize}[nosep]
\item[1.] We study workload balance among drivers over an extended planning horizon. With this viewpoint, we show that workload equity can be achieved with very limited impact on economic efficiency.
\item[2.] We formally introduce the MVRPB and a simple and efficient two-phase approach for its solution. This approach also permits obtaining bounds on the best possible workload equity.
\item[3.] Through extensive numerical experiments, we demonstrate the performance of the proposed solution approach. Moreover, we measure (i) the gap between the equity level achieved by the two-phase method and perfect equity and (ii) the benefits of considering a longer planning horizon.
\end{itemize}

The remainder of this paper is organized as follows. Section~\ref{LR} reviews related works for single and multi-period VRPs. Section~\ref{PS} defines the MVRPB, and Section~\ref{sec:math_form} describes the proposed two-phase solution approach. Section~\ref{numeric anal} presents a computational study based on test instances derived from well-known capacitated VRP benchmark instances. Finally, Section~\ref{conclusion} concludes the paper and discusses research perspectives.

\section{Related Studies}
\label{LR} 

\noindent
\textbf{Single-period planning.}
The classical capacitated VRP (CVRP) seeks a set of routes starting and ending at a central depot and visiting a given set of clients \citep{Vidal2019ACG}. Each client is characterized by a demand quantity and must be serviced in a single visit. The total demand of the clients over each route should not exceed the vehicle's capacity (considered to be identical for all vehicles). The objective of the problem is to minimize the total traveled distance. The CVRP is NP-hard as it generalizes the traveling salesman problem (TSP). Consequently, no optimal (i.e., exact) polynomial-time algorithm is known for its solution. Moreover, the best existing exact approaches can only solve medium-sized problems counting a few hundred customers in less than a few hours. Accordingly, heuristic approaches are widely used in practical settings \citep{articleRibeiro, VIDAL20131}. 

 \cite{BOWERMAN1995107} were among the first to introduce an equity objective in the context of a school bus routing problem. They used a variance measure to balance the length of the drivers' routes and the number of students transported. In addition, they considered two additional criteria within a multi-objective framework. Since that work, numerous papers have considered workload balance objectives in single-period VRPs with different equity measures (see, e.g., \citealt{GOLDEN1997445,lee1999study,articleJozefowiez2009,articleLopez,oyola2014grasp,bertazzi2015min, articleGalindres, articleVa, articleVega-Mejia, zhang2019gmma, LEHUEDE2020129,  articleLondono}). The \textsc{Min-Max} measure and the difference between the maximum and minimum workloads (called \textsc{Range} by \citealt{Matl_2018}) are most commonly used. When considering distance as the workload metric, \textsc{Min-Max} minimizes the longest route in a solution (see, e.g., \citealt{articleLopez}), whereas \textsc{Range} minimizes the difference between the longest and shortest route (see, e.g., \citealt{articleLondono}). We refer the reader to \cite{HALVORSENWEARE2016451}, \cite{articleLozano}, and \cite{Matl_2018} for a comprehensive list of equity measures.

No solution simultaneously optimizes cost and equity in most situations, especially when considering equity within a single period. Consequently, equity is often considered within multi-objective approaches, in which trade-off solutions have to be found between equity and other objectives such as distance. Many studies along this line refer to the resulting problem as the VRP with ``route balancing'' \citep{articleJozefowiez2007, articleJozefowiez2009}. Some studies have examined how equity objectives affect the shape of the Pareto front in a bi-objective context \citep{articleRaul, HALVORSENWEARE2016451, articleLozano, articleSchwarze, zhang2019gmma}. \cite{Matl_2018, articleMatl2019} extensively discussed six equity measures based on eight desirable properties. They claim that no measure satisfies every property or is strictly better than the others for all relevant properties. They also show that a monotonic equity measure such as \textsc{Min-Max}, which is based on a variable-sum metric such as distance, avoids non-TSP-optimal solutions and inconsistent solutions and therefore should be an objective of choice. A solution is non-TSP optimal if at least one route can be rearranged while improving its distance. Inconsistency occurs when a given solution is preferred over another solution even though all of its routes are longer. In both cases, the distance of one or more routes has been artificially increased to improve the equity objective.

\noindent
\textbf{Multi-period planning.}
Multi-period VRPs are defined over a time horizon of several periods and generally aim to minimize the total routing cost over all periods. In the periodic VRP (PVRP), each client is characterized by a visit frequency representing how many visits are requested over the planning horizon and a list of patterns representing acceptable visit-day combinations. Solving this problem requires selecting a visit pattern for each client and generating the routes for each period. As a consequence, the decisions made at each period become interdependent \citep{CoeneSofieperiodic}.
\cite{mourgaya2006periodic} presented several PVRP variants and classified them based on their objectives, constraints, and solution methods. Equity is one of the objectives discussed in that study.

Some studies considered multi-period VRPs and measured equity within each period. Papers in this category generally arise from real-life case studies and involve different workload equity objectives. \cite{articleBlakeley} studied the problem of an elevator-maintenance company that assigns technicians to clients. The objective was to minimize a weighted sum of travel time, overtime, and unbalanced workload in each period. To reach equity, the authors designed heuristics that optimize a \textsc{Range} measure on the travel times. \cite{10.5555/1670696.1670704} studied a multi-period VRP to balance the workload of meter readers over a month for a utility company. To achieve a good balance, they set a lower and upper bound on the number of clients and the length of each daily route. They also constrain the deviation of a client’s bill from one month to the next. Their three-stage methodology combines heuristics and integer programming. \cite{GULCZYNSKI2011648} presented a PVRP in which an equal workload is sought within each period. The overall objective is a weighted sum of distance and \textsc{Range} measure over the number of clients served. The problem was solved with an integer-programming-based heuristic. 

\cite{mourgaya2006periodic} designed a hierarchical heuristic to solve an industrial PVRP application with 16\,658 visits over a time horizon of 20 days. At a tactical level, the method allocates each client to a combination of days and vehicles. The objective considered at this stage involves minimizing the maximal workload over the days and vehicles subject to additional geographical restrictions. At the operational level, the method minimizes the traveled distance.

\cite{MOURGAYA20071028} then studied a tactical planning problem that required choosing visit days for clients subject to service level constraints. Once the visit days were selected, the clients were assigned to vehicles to achieve equity among drivers. At this stage, equity was modeled by a constraint that bounds the workload (served demand) associated with each cluster of clients. In another work, \cite{articleLinfati} developed a two-phase heuristic to balance the number of medication deliveries to patients among the drivers for each day of the planning horizon. They considered a \textsc{Range} measure within a weighted-sum objective that also accounts for extra hours, extra capacity, as well as daily and client clustering costs. Finally, \cite{inbook} introduced a PVRP with different possible visit frequencies for each client. They relied on objectives that minimize the total route cost and the number of stops to the same client, and included a workload balance component between vehicles for each period. This component of the objective minimizes the maximum service time over the weeks of the planning horizon and the days of the week.
Furthermore, they imposed a maximum route duration for each vehicle. A decomposition approach was used first to assign clients to weeks, and then to assign them to a day within the selected week. A variant of the classical VRP was finally solved for each day using a three-phase adaptive large-neighborhood search. This algorithm was applied in a case study for a hygiene service company performing more than 69\,951 visits for 6000 clients over 12 weeks, leading to significant practical savings.

\cite{6548751} studied a PVRP for home health care in which three types of patients require services. The primary objective was to reach equity by minimizing the maximum route time for all vehicle routes over the week. The solution method combined tabu search and local search. \cite{8968736} balanced the workload among routes for a periodic home health care assignment problem by partitioning the service area into regions. Finally, \cite{schonberger2016multi} set upper bounds on route duration to ensure daily workload balance in a multi-period VRP.
 
 \cite{articleRibeiro} was, to our knowledge, one of the very few studies considering a multi-period VRP in which equity is evaluated over multiple periods. The problem is solved with a weighted sum objective considering cost, equity, and market share. The equity objective minimizes the standard deviation of the workloads over many periods, where the workload corresponds to the demand served by each driver. Small instances were solved with a commercial solver. \cite{articleHuang} studied a PVRP for which the objective function minimizes the total workload of all drivers, and the maximum workload difference between two drivers cannot exceed a threshold. The equity metric, in this case, corresponds to the sum of travel and service times.
 The results show that workload equity among drivers can be achieved at a reasonable cost. Finally, \cite{articleMancini} recently considered workload balance and service consistency in a collaborative multi-period VRP, where the number of clients assigned to a carrier over the planning horizon is constrained. They also examined how workload balance and service consistency impact the total solution cost.

As seen in this review, very few studies have focused on multi-period VRPs that account for equity, and even fewer studies have considered equity over multiple periods. The present study fills this critical methodological gap. It proposes a simple solution approach for an equitable vehicle routing problem defined over multiple periods, and analyses the resulting equity depending on the planning horizon.

\section{Problem Statement}
\label{PS}

The MVRPB can be formally defined on a complete graph $G$ = ($V$, $E$), where $V$ is the set of vertices and $E$ is the set of edges.
Let $V = \{0\} \cup C$, where $C$ is a set of vertices representing clients and vertex 0 is the depot.
The cost $d_{ij}$ corresponds to the length of edge $(i,j) \in E$ representing a direct trip from $i$ to $j$.
For simplicity, we assume that all distances $d_{ij}$ are integer. Each client $i \in C$ is characterized by a list of visit days (i.e., periods) within a planning horizon of $T$ periods and by a demand quantity $q^t_{i}$ on each of these days (i.e., demand quantities can differ between days). Finally, $m$ drivers are available through the planning horizon to perform the visits, and we are given an unlimited fleet of homogeneous vehicles with capacity $Q$ located at the depot.

Solving the MVRPB amounts to finding routes for each day in such a way that (i) each client is visited on each requested day, (ii) each route for a given day is assigned to a single driver, (iii) no driver operates more than one route in a day, and (iv) no route exceeds the vehicle capacity. Note that drivers may not necessarily work during each day of the planning horizon. The objectives are to optimize distance and balance the drivers' workloads over the planning horizon.

\cite{Matl_2018} discussed different \emph{equity metrics} (e.g., travel time, distance, demand quantity served, number of clients) and \emph{equity measures} (e.g., the maximum workload of a driver, the difference between the smallest and largest workloads, and the standard deviation of workloads). In this work, we focus on distance-based workload equity among drivers. Therefore, the workload of a driver is the total distance traveled by the driver over the planning horizon. We use the \textsc{Min-Max} equity measure, which aims to achieve equity by minimizing the maximum total distance of the routes traveled by any driver over the planning horizon.

\section{Solution Approach}
\label{sec:math_form}

To solve the MVRPB, we design a two-phase solution approach. As seen in the remainder of this section, our approach follows a hierarchical objective: it first guarantees a routing solution with minimum cost (i.e., with minimum total distance over the planning horizon) and then maximizes workload balance. This permits us to evaluate the extent to which workload equity can be ensured over a longer planning horizon without sacrificing economic efficiency. Moreover, as seen in the following, the optimal routes found in the first phase will permit us to calculate a bound on the best possible workload equity for any solution.

Our method unfolds in two stages. Firstly, it solves a CVRP for each period considering only the deliveries of this period and minimizing cost (total distance). Then, it solves an allocation problem to assign routes to drivers on each period, intending to optimize workload balance. The remainder of this section details the techniques designed to perform each step efficiently.

\subsection{First Phase -- Distance Optimization}

The first phase consists of solving one CVRP per period $t \in T$ to obtain high-quality routing plans minimizing distance. We rely on mixed-integer linear programming (MILP) techniques to solve these problems. We first present a compact mathematical formulation of the problem and then discuss its solution by branch-and-price.

\begin{table}[htbp]
\caption{Notations used in the mathematical programs} \label{tab:notations}
\vspace*{-0.25cm}
\centering
\scalebox{0.85}
{
\begin{tabular}{l @{\hspace*{1cm}} l @{\hspace*{1cm}} l}
\toprule
\textbf{Sets}  & $C$ & Set of clients\\ 
& $E$ & Set of edges \\
& $V$ & Set of nodes, $V = C \cup \{0\}$\\
& $C_t$ & Set of clients in period $t \in T$\\
& $E_t$ & Set of edges in period $t \in T$ \\ 
& $V_t$ & Set of vertices in period $t \in T$, $V_t = C_t \cup \{0\}$\\
& $R_t$ & Set of routes in period $t \in T$ \\ 
\midrule
\textbf{Parameters} & $T$ & Number of time periods \\ 
 & $m$ & Number of drivers (i.e., maximum number of vehicle routes in each period) \\
 & $Q$ & Capacity of each vehicle \\
& $q_{i}^{t}$ & Demand of client $i \in C_{t}$ in period $ t \in T$ \\ 
& $d_{ij}$ & Distance of edge $(i,j) \in E$ \\ 
\midrule
\textbf{Variables} & $y_{ir}$ & Binary variable equal to 1 if client $i$ is on route $r$, $i \in C$, $r \in R_t$. \\ 
& $x_{ijr}$ & Binary variable equal to 1 if edge $(i,j) \in E$ is on route $r$, $i,j \in V$, $r \in R_t$. \\ 
&& $r \in R_t$, $t \in T$\\
\bottomrule
\end{tabular}
}
\end{table}

Table \ref{tab:notations} summarizes the notations used in the mathematical models.
For each period~$t$, the resulting CVRP$_t$ subproblem can be mathematically formulated using a three-index undirected formulation, as in \cite{articleBaldacci}. In this formulation, the cut set for any $S \subseteq C_t$ is defined as $\delta(S)= \{(i, j) \in E_t : i \in S$, $j \notin S$ or $i \notin S $, $j \in S \}$:
\begin{align}
(CVRP_t) \hspace*{0.7cm} \min & \sum_{r=1}^m \sum_{(i,j) \in E_t} d_{ij} x_{ijr} & \label{eq:base_obj0}\\
\mbox{s.t.} \ & \sum_{r=1}^m y_{ir} = 1 & i \in C_t & \label{eq:base_cust assignment0}\\
&\sum_{i \in C_t} q_{i}^{t} y_{ir}\leq Q
& r\in\{1,\dots,m\} & \label{eq:base_capacity0}\\
&\sum_{(i,j)\in \delta(\{i\})} x_{ijr} = 2 y_{ir} & i \in C_t, r\in\{1,\dots,m\} & \label{eq:base_folow0} \\
&\sum_{(i,j)\in \delta(S)}x_{ijr} \geq y_{ir} &S \subseteq C_t \, : \, \vert S \vert \geq 2, i \in S, r\in\{1,\dots,m\} &\label{eq:subtour0}\\ 
& x_{ijr} \in \{0, 1\} & (i,j) \in E_t \backslash \{(0,j):j \in C_t\} , r\in\{1,\dots,m\}& \label{eq:binary0} \\
& x_{0jr} \in \{0, 1, 2\} & \ j \in C_t , r\in\{1,\dots,m\} & \label{eq:depot_binary0} \\
& y_{ir} \in \{0, 1\} & i \in C_t, r\in\{1,\dots,m\} &\label{eq:binary00}
\end{align}

Objective (\ref{eq:base_obj0}) minimizes the total routing cost. Constraints (\ref{eq:base_cust assignment0}) force each client to be served by exactly one route. Constraints (\ref{eq:base_capacity0}) ensure that capacity constraints are respected. Constraints (\ref{eq:base_folow0}) are flow-conservation constraints for the routes that also link the $x$ and $y$ variables. Constraints (\ref{eq:subtour0}) are the sub-tour-elimination constraints. These constraints guarantee that each route contains the depot.

Previous studies on integer programming approaches for the CVRP indicate that directly solving Model~(\ref{eq:base_obj0}--\ref{eq:binary00}) is ineffective for medium instances with more than a few dozen clients and that combining cuts and column generation generally provides better results. Accordingly, we exploit the VRPSolver framework \citep{vrpsolver} for an efficient solution to this problem. This solver provides a generic branch-and-cut-and-price (BCP) framework for many classes of MILPs, including, among others, the considered problem setting. It has achieved a competitive or superior performance on standard test instances when compared with specialized VRP algorithms and is currently available at \url{https://vrpsolver.math.u-bordeaux.fr/}.

This algorithmic framework relies on the solution of successive pricing subproblems that take the form of resource-constrained shortest paths (RCSPs) on a path-generator graph (VRPGraph). A bidirectional-labeling dynamic programming algorithm is then used to solve the RCSPs. 
VRPSolver relies on the concept of packing sets to generalize well-known cuts. Essentially, a packing set is a subset of arcs such that at most one arc from the given subset appears in the paths that are part of an optimal solution. Packing sets are defined in accordance with the application considered and associated model. In our specific case, the packing sets represent limited memory rank-1 cuts (a generalization of the subset row cuts -- \citealt{articleJepsen}) and rounded capacity cuts \citep{Laporte1983ABA}. The branching rule in VRPSolver is based on accumulated resource consumption and, if needed to enforce integrality, on a generalization of the branching rule of \citet{ryan1981rn}.

Finally, the performance of VRPSolver depends on the availability of a good initial upper bound (UB) to limit the search. To find such an initial solution and bound, we use the hybrid genetic search (HGS -- \citealt{vidal2012hybrid, HGS}), a state-of-the-art metaheuristic for the CVRP. HGS could be used as a stand-alone approach for the first phase in time-critical applications, or if the scale of the problems becomes too large for an exact solution. Still, we opted to additionally rely on the exact algorithm in this phase, as this will allow us to obtain lower bounds on the best achievable workload equity (see next section).

\subsection{Second Phase -- Equitable Workload Allocation}
\label{second step model}

The second phase of the algorithm takes as input the solution $R^*$ found in the previous phase, represented as the set of optimal routes $R^*_t$ for each period. It seeks to achieve a fair distribution of these routes among $m$ identical drivers. In this stage, we use a \textsc{Min-Max} objective to minimize the maximum workload of any driver (i.e., the maximum total distance driven by a driver over all periods).

\subsubsection{MILP formulation and bounds}
\label{LB}

This allocation problem can also be mathematically formulated as a MILP. Let $d_r$ represent the distance driven on each route $r \in R^*_t$ found in the first phase, and let $z^{t}_{rk}$ be a binary variable equal to $1$ if route $r$ is assigned to driver $k$ in time period $t$. Finally, let $\Delta$ be a continuous variable capturing the maximum distance for a driver over the planning horizon. The best possible workload equity for the considered set of routes can be found by solving the following model:
\begin{align}
&\min \ \Delta &\label{eq:equity obj1} \\
\mbox{s.t.}\ & \sum_{t=1}^T \sum_{r \in R^*_t} d_{r} z^{t}_{rk} \leq \Delta & k \in \{1,\dots,m\} &\label{eq: driver workload1}\\
& \sum_{r \in R^*_t} z^{t}_{rk} \leq 1 & t \in \{1,\dots,T\}, k\in \{1,\dots,m\}  & \label{eq:assignment of drivers to routes1} \\
&\sum_{k=1}^m z^{t}_{rk} = 1 & t \in \{1,\dots,T\}, r\in R^*_t & \label{eq:assignment of routes to drivers1}\\
& z^{t}_{rk} \in \{0, 1\} &  t \in \{1,\dots,T\}, r\in R^*_t, k \in \{1,\dots,m\} \label{binary1}&\\
& \Delta \in \mathbb{R}^+. \label{free1}&
\end{align}

Objective~(\ref{eq:equity obj1}) and~Constraints (\ref{eq: driver workload1}) model the minimization of the maximum workload over all drivers. Constraints (\ref{eq:assignment of drivers to routes1}) ensure that each driver serves at most one route in each period, whereas Constraints (\ref{eq:assignment of routes to drivers1}) ensure that each route is assigned to exactly one driver in each period. Finally, Constraints (\ref{binary1}) and (\ref{free1}) define the domain of the decision variables. Note that the inequality in Constraints (\ref{eq:assignment of drivers to routes1}) can be replaced by an equality if the number of drivers matches the number of routes found in the first phase in any given period~$t$. The resulting formulation can be viewed as a variant of the bin packing problem (BPP) with conflicts~\citep{Capua2018}.

Let $\Delta_\textsc{opt}(R)$ be the optimal workload produced by solving the allocation problem for a given first-phase routing solution $R$, and let $D(R)$ be the total distance of all routes of a solution $R$ over all periods. The following bounds are valid:

\begin{property}
\label{P1}
The best possible workload allocation for a routing solution $R$ is such that
\begin{equation}
\Delta_\textsc{opt}(R) \geq \left\lceil \frac{D(R)}{m} \right\rceil.
\end{equation}
\end{property}

\noindent
\textbf{Proof.}
The proof directly derives from the formulation of the assignment problem. Summing Constraint~(\ref{eq: driver workload1}) over $k \in \{1,\dots,m\}$ gives:
\begin{equation}
\sum_{k=1}^m \sum_{t=1}^T \sum_{r \in R_t} d_{r} z^{t}_{rk} \leq m \Delta_\textsc{opt}(R).
\end{equation}
Next, using Equation~(\ref{eq:assignment of routes to drivers1}) leads to:
\begin{equation}
\sum_{t=1}^T \sum_{r \in R_t} d_{r}  \leq m \Delta_\textsc{opt}(R) \Longrightarrow \Delta_\textsc{opt}(R) \geq \frac{D(R)}{m}
\end{equation}
Finally, since the distances $d_{ij}$ are integer, then $\Delta_\textsc{opt}(R)$ is also an integer, and the right-hand side of the inequality can be rounded up, giving the announced result.
 
\begin{property}
\label{P2}
Let $\Delta_\textsc{opt}$ be the best possible workload equity achievable in any solution of the MVRPB (including first-stage solutions that are not optimal in terms of distance), then:
\begin{equation}
\label{eq:optbound}
\Delta_\textsc{opt} \geq \left\lceil \frac{D(R^*)}{m} \right\rceil.
\end{equation}
\end{property}

\noindent
\textbf{Proof.}
As a consequence of Property~\ref{P1}, the best possible workload equity over all possible routing solutions $R$ satisfies the following relation:
\begin{equation}
\Delta_\textsc{opt} = \min_{R} \Delta_\textsc{opt}(R) \geq \min_{R} \left\lceil \frac{D(R)}{m} \right\rceil = \left\lceil  \frac{ \min_{R} D(R)}{m}  \right\rceil = \left\lceil \frac{D(R^*)}{m} \right\rceil.
\end{equation}

This relation gives us a lower bound $\text{LB} = \left\lceil D(R^*)/m \right\rceil$, which permits us to evaluate how far our MVRPB solutions are from the best possible workload equity, calculated by assuming that the total amount of workload from distance-optimal routing solutions is evenly distributed among drivers. It is important to remark that $\Delta_\textsc{opt}$ can be smaller than $\Delta_\textsc{opt}(R^*)$ since the best workload balance over multiple periods can involve routes that do not belong to any optimal CVRP solution of a given period. In contrast, the bound announced in Property 2 holds for any MVRPB solution.
 
 \subsubsection{Set-partitioning reformulation and solution approach}
 \label{UB}
 
A direct solution of Formulation~(\ref{eq:equity obj1}--\ref{free1}) using standard MILP solvers such as CPLEX is ineffective. This is partly due to symmetry, given that all drivers are considered identical, and renumbering them produces equivalent solutions. One way to circumvent this issue is to solve this problem as a sequence of set-partitioning feasibility problems for different values of $\Delta$ as defined in Model~(\ref{SP1}--\ref{SP3}).
\begin{align}
& \sum_{\sigma \in \Omega_\Delta} \lambda_\sigma = m \label{SP1}\\
& \sum_{\sigma \in \Omega_\Delta} a_{\sigma r} \lambda_\sigma = 1 & t \in \{1,\dots,T\}, r \in R^*_t \label{SP2}\\
& \lambda_\sigma \in \{0, 1\} &  \sigma \in \Omega_\Delta \label{SP3}
\end{align}

In this formulation, each element $\sigma \in \Omega_\Delta$ represents an admissible combination of routes (i.e., a schedule) that a driver can operate over the planning horizon without exceeding a workload of~$\Delta$. Each constant $a_{\sigma r}$ takes value $1$ if and only if $r$ belongs to $\sigma$. Each binary variable $\lambda_\sigma$ takes value~$1$ if this schedule is selected for one driver. Constraints~(\ref{SP1}) ensure that work schedules are created for the $m$ drivers, and Constraints~(\ref{SP2}) ensure that each route appears exactly in one schedule. 

\noindent
\textbf{Binary search strategy.}
Finding a feasible solution of Model~(\ref{SP1}--\ref{SP3}) for a given $\Delta$ means that there exists a feasible allocation of the routes to drivers in such a way that the maximum workload over the planning horizon does not exceed~$\Delta$. Therefore, the optimal workload balance is such that $\Delta_\textsc{opt}(R^*) \leq \Delta$. In contrast, proving that this model is infeasible would imply that $\Delta_\textsc{opt}(R^*) > \Delta$.

With this, we develop a strategy based on a binary search to locate $\Delta_\textsc{opt}(R^*)$. Initially, we start with $\Delta= LB = \left\lceil \frac{D(R^*)}{m} \right\rceil$ and solve Model~(\ref{SP1}--\ref{SP3}). If this model is feasible, then we have attained the best possible workload equity. Otherwise, this means there is no solution with a workload balance of $\Delta = \left\lceil \frac{D(R^*)}{m} \right\rceil$. In this case, we use a \textbf{construction heuristic} to find an initial feasible solution and therefore an upper bound (UB) value for $\Delta$. This heuristic consists in ordering the items (i.e., all the routes over the planning horizon) by decreasing workload (distance) and then following this order to insert them iteratively into a compatible bin (i.e., a driver that does not operate a route on that day) that has the current smallest workload.
At this point, we have a range for the optimum (integer) value and locate it by binary search, using Model~(\ref{SP1}--\ref{SP3}) to determine feasibility at each step. This process stops when the model for $\Delta_\textsc{opt}(R^*)-1$ is infeasible and the model for $\Delta_\textsc{opt}(R^*)$ is feasible.

\noindent
\textbf{Solution of each subproblem}.
Model~(\ref{SP1}--\ref{SP3}) contains an exponential number of variables~$\sigma \in \Omega_\Delta$, therefore a direct solution approach is impractical. To solve this problem, we rely once again on the branch-and-price framework provided by VRPSolver. \cite{pessoa:hal-02986956} provides adaptations of VRPSolver to the classical BPP and other variants such as vector packing, variable-sized BPP, and variable-sized BPP with optional items. To solve our problem with VRPSolver, we essentially need to redefine the path-generator graph (VRPGraph).

\begin{figure}[htbp]
\centering
\begin{tikzpicture}[-latex ,auto, node distance={30mm}, thick, main/.style = {draw, circle}] 
\node[main] (1) {$v_0$}; 
\node[main] (2) [right of=1] {$v_1$};
\node[main] (3) [right of=2] {$v_2$};
\node[main] (4) [right of=3] {$v_3$};
\node[main] (5) [right of=4] {$v_{I-1}$};
\node[main] (6) [right of=5] {$v_I$};
\path (1) edge [bend left =45] node[above] {$a_{1+}$ }(2);
\path (1) edge [bend left =-45] node[above] {$a_{1-}$}(2);
\path (2) edge [bend left =45] node[above] {$a_{2+}$ }(3);
\path (2) edge [bend left =-45] node[above] { $a_{2-}$}(3);
\path (3) edge [bend left =45] node[above] {$a_{3+}$ }(4);
\path (3) edge [bend left =-45] node[above] { $a_{3-}$}(4);
\node[main] (7) [right of=4,draw=none,node distance=1.4cm] {$\ldots$};
\path (5) edge [bend left =45] node[above] {$a_{I+}$}(6);
\path (5) edge [bend left =-45] node[above] {$a_{I-}$}(6);
\end{tikzpicture}
\vspace{-0.2cm}
\caption{Path-generator graph for the BPP.}
\label{paper graph}
\end{figure}
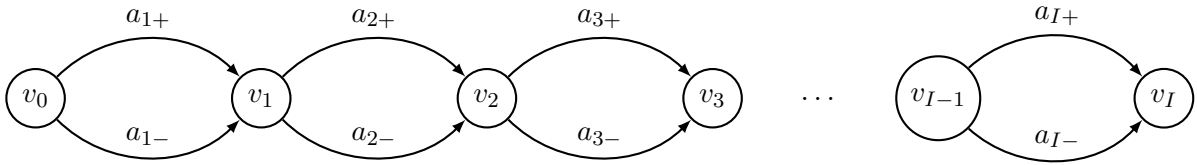

Figure~\ref{paper graph} represents the VRPGraph for the classical BPP, assuming that $I$ items need to be packed. Each node in this graph corresponds to one item, except node $v_0$, representing a starting point. Item~$i$ of weight $w_i$ is loaded in the bin each time we use arc $a_{i+}$ along the path from the start node to the end node (conversely, item $i$ is not loaded in the bin if arc $a_{1-}$ is used). Next, each path generated in this graph that does not exceed the bin capacity defines a new packing (column) for the column-generation algorithm. Capacity is the only resource consumed in VRPGraph, and a packing set is made of a subset of arcs with nonzero consumption of the resource.

\begin{figure}[htbp]
\centering
\begin{tikzpicture}[-latex ,auto, node distance={30mm}, thick, main/.style = {draw, circle}] 
\node[main] (1) {$P_0$}; 
\node[main] (2) [right of=1] {$P_1$};
\node[main] (3) [right of=2] {$P_2$};
\node[main] (4) [right of=3] {$P_3$};
\node[main] (5) [right of=4] {$P_4$};
\node[main] (6) [right of=5] {$P_5$};
\path (1) edge [bend left =15] node[above = -0.1cm] {\footnotesize $a_{11+}$ }(2);
\path (1) edge [bend left =45] node[above=-0.1cm] {\footnotesize $a_{12+}$} (2);
\path (1) edge [bend left =60 ,draw=none] node[above] {\footnotesize$\ldots$ }(2);
\path (1) edge [bend left =89] node[above] {\footnotesize $a_{1n_1+}$}(2);
\path (1) edge [bend left =-75] node[above] {\footnotesize $a_{1-}$}(2);
\path (2) edge [bend left =15] node[above = -0.1cm] {\footnotesize $a_{21+}$}(3);
\path (2) edge [bend left =45] node[above=-0.1cm] {\footnotesize $a_{22+}$}(3);
\path (2) edge [bend left =60 ,draw=none] node[above] {\footnotesize$\ldots$ }(3);
\path (2) edge [bend left =89] node[above] {\footnotesize $a_{2n_2+}$}(3);
\path (2) edge [bend left =-75] node[above] {\footnotesize $a_{2-}$}(3);
\path (3) edge [bend left =15] node[above = -0.1 cm] {\footnotesize $a_{31+}$}(4);
\path (3) edge [bend left =45] node[above=-0.1cm] {\footnotesize $a_{32+}$}(4);
\path (3) edge [bend left =60 ,draw=none] node[above] {\footnotesize$\ldots$ }(4);
\path (3) edge [bend left =89] node[above] {\footnotesize $a_{3n_3+}$}(4);
\path (3) edge [bend left =-75] node[above] {\footnotesize $a_{3-}$}(4);
\path (4) edge [bend left =15] node[above = -0.1cm] {\footnotesize $a_{41+}$}(5);
\path (4) edge [bend left =45] node[above =-0.1cm] {\footnotesize$a_{42+}$}(5);
\path (4) edge [bend left =60 ,draw=none] node[above] {\footnotesize$\ldots$ }(5);
\path (4) edge [bend left =89] node[above] {\footnotesize $a_{4n_4+}$}(5);
\path (4) edge [bend left =-75] node[above] {\footnotesize $a_{4-}$}(5);
\path (5) edge [bend left =15] node[above=-0.1cm] {\footnotesize $a_{51+}$}(6);
\path (5) edge [bend left =45] node[above=-0.1cm] {\footnotesize $a_{52+}$}(6);
\path (5) edge [bend left =60 ,draw=none] node[above] {\footnotesize$\ldots$ }(6);
\path (5) edge [bend left =89] node[above] {\footnotesize $a_{5n_5+}$}(6);
\path (5) edge [bend left =-75] node[above] {\footnotesize $a_{5-}$}(6);
\end{tikzpicture} 
\vspace{-0.2cm}
\caption{Path-generator graph for the second-stage model.}
\label{my graph}
\end{figure}
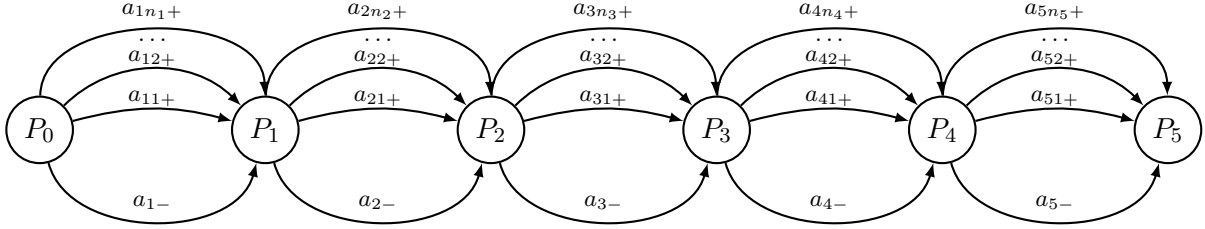

Figure \ref{my graph} provides an adapted path-generator graph for Problem~(\ref{SP1}--\ref{SP3}). In this graph, nodes correspond to periods, except the first node $P_0$, which represents a starting point. An arc of type $a_{ij+}$ then goes from $P_{i-1}$ to $P_i$ for each possible route $j \in \ \{1, 2, 3, \dots , n_i\}$ in period $P_i$, where $n_i$ is the number of routes in period~$P_i$. If one of these parallel arcs is used, the corresponding route is assigned to the driver and contributes to the driver's total workload (total distance limited to $\Delta$, which stands as the bin's capacity). An arc of type $a_{i-}$ is also available to represent the possibility of assigning no route to the driver in period~$P_i$. With these conventions, a path from the start node to the end node that does not exceed the workload limit $\Delta$ corresponds to a feasible assignment of routes to a driver over the planning horizon. 

We use the standard parameter setting of VRPSolver with just a minor modification to its diving heuristic. VRPSolver typically uses a diving heuristic before branching to improve the primal solution \citep{sadykov2019primal}, but only at the root node. In contrast, our implementation allows strong diving at each node to quickly locate feasible solutions. If such a solution is found, then the solver can be immediately stopped since the model is known to be feasible.

 \section{Computational Study}
 \label{numeric anal}
 
The goal of our experimental study is twofold: 
(i) evaluating the performance of the proposed solution approach and especially the computational effort needed for each of its steps, and 
(ii) measuring to which extent workload equity can be achieved over multiple periods without sacrificing economic efficiency, by simply allocating cost-optimal routes to drivers in an equitable fashion as done in our two-phases approach.

Our experiments are conducted on a 2.4\,GHz Intel Gold 6148 Skylake processor with 8\,GB of RAM. The VRPSolver interface is implemented in Julia v1.4.2 with JuMP v0.18. VRPSolver uses BaPCod, a C++ library for implementing a generic BCP, and CPLEX 12.8 to solve the linear and mixed-integer linear programs. 
All experiments have been conducted on a single thread.

\subsection{Test Instances}
\label{test case instances}
 
To construct test instances for the MVRPB, we rely on a subset of the CVRP instances of \cite{articleUchoa}, as they include diverse characteristics: distribution and number of clients, depots locations, and average route length. The complete set contains 100 Euclidean instances with 100 to 1000 clients. The distances are rounded to the nearest integer as in the original instances.

Since the MVRPB is defined in a multi-period context, we had to modify the original instances. Therefore, we selected the instances X-n200-k36, X-n204-k19, X-n209-k16, X-n214-k11, X-n219-k73, X-n223-k34, X-n228-k23, X-n233-k16, X-n237-k14, and X-n242-k48 including between 199 and 241 clients from \cite{articleUchoa}. We generated three different 10-period MVRPB configurations for each of these ten instances by randomly selecting 50, 75, or 100 clients in each period. Finally, for each period and client with demand $d$ in the original instance, we randomly selected a new demand realization from a uniform integer distribution in $\{\lceil 0.5\times d \rceil\, \dots, \lceil 1.5 \times d \rceil \}$. This way, clients can have different demands at different periods. We repeated this generation (customers and demands selection) ten times for each configuration, leading to $10 \times 3 \times 10 = 300$ MVRPB instances defined over $10$ periods.
Finally, to obtain instances with fewer periods $T \in \{2, 3, 5, 7\}$, we retained the first $T$ periods of each 10-period instance. The number of drivers in each 10-period instance was set to the maximum number of routes from optimal CVRP solutions over the ten periods. Consequently, some drivers may be idle for a given period. The number of drivers in the 2-, 3-, 5-, and 7-period instances is kept identical to the number of drivers in the corresponding 10-period instance.

\subsection{Computational Performance}
\label{P and S}

We first evaluate the performance of each of the two steps of the proposed approach: the computational effort needed to find optimal CVRP solutions in each period, and the effort to find an equitable workload allocation in the second step. We refer to these steps as (1) route optimization and (2) multi-period workload balancing.\\

\noindent
\textbf{Route optimization.}
Tables~\ref{tab1}~to~\ref{tab3} report the performance of the route-optimization step for the 10-period instances with 50, 75, and 100 clients in each period. For brevity, the results are presented in aggregated form, with one line for each original instance of \cite{articleUchoa}, by averaging over the $10$ corresponding MVRPB instances and $10$ periods. From left to right, the columns report the names of the associated original instances, the average traveled distance per period in the solutions found by HGS and VRPSolver, the average computational time of these two methods, and finally, the number $\#k$ of drivers. 

\begin{table}[htbp]
 \fontsize{10}{12}\selectfont
\caption{Performance of the route-optimization step, for MVRPB instances with 50 clients per period} \label{tab1}
 \centering
 \setlength\tabcolsep{0.15cm}
  \scalebox{0.95}
 {
  \begin{tabular}{c@{\hspace*{1cm}}rr@{\hspace*{1cm}}rr@{\hspace*{1cm}}r}
 \hline
         \multirow{2}{*}{Instance} & \multicolumn{2}{c}{Distance} & \multicolumn{2}{c}{T(s)} &\multirow{2}{*}{\#k}  \\ 
           \cmidrule(lr){2-3}\cmidrule(lr){4-5}
       & HGS & VRPSolver & HGS & VRPSolver & \\ \hline
       
        X-n200-k36 & 16284.45 & 16281.48 & 15.8 & 614.1 & 10.0 \\ 
        X-n204-k19 & 7200.85 & 7200.85 & 14.0 & 7.8 & 5.3 \\ 
        X-n209-k16 & 9699.69 & 9699.69 & 14.2 & 14.2 & 4.8 \\ 
        X-n214-k11 & 3782.88 & 3782.88 & 14.2 & 52.1 & 3.4 \\ 
        X-n219-k73 & 28562.15 & 28562.15 & 13.4 & 2.3 & 17.0 \\ 
        X-n223-k34 & 11464.69 & 11464.17 & 13.1 & 7.3 & 9.5 \\ 
        X-n228-k23 & 7652.85 & 7652.72 & 14.1 & 25.5 & 7.3 \\ 
        X-n233-k16 & 6704.64 & 6704.64 & 12.8 & 12.9 & 4.5 \\ 
        X-n237-k14 & 7980.45 & 7980.45 & 13.1 & 31.4 & 3.0 \\ 
        X-n242-k48 & 19970.99 & 19970.05 & 15.0 & 7.6 & 12.1 \\ \hline
 \end{tabular}
 }
\end{table}

\begin{table}[htbp]
 \fontsize{10}{12}\selectfont
\caption{Performance of the route-optimization step, for MVRPB instances with 75 clients per period} \label{tab2}
 \centering
 \setlength\tabcolsep{0.15cm}
  \scalebox{0.95}
 {
  \begin{tabular}{c@{\hspace*{1cm}}rr@{\hspace*{1cm}}rr@{\hspace*{1cm}}r}
 \hline
         \multirow{2}{*}{Instance} & \multicolumn{2}{c}{Distance} & \multicolumn{2}{c}{T(s)} &\multirow{2}{*}{\#k}  \\ 
           \cmidrule(lr){2-3}\cmidrule(lr){4-5}
       & HGS & VRPSolver & HGS & VRPSolver & \\ \hline
       X-n200-k36 & 23388.65 & 23379.47 & 28.7 & 185.5 & 14.7 \\ 
        X-n204-k19 & 9294.97 & 9294.97 & 22.0 & 31.6 & 7.9 \\ 
        X-n209-k16 & 13272.78 & 13272.68 & 24.5 & 60.1 & 6.6 \\ 
        X-n214-k11 & 4881.22 & 4881.2 & 25.2 & 268.2 & 5.0 \\ 
        X-n219-k73 & 41826.76 & 41826.76 & 20.3 & 2.8 & 25.0 \\ 
        X-n223-k34 & 16035.39 & 16034.47 & 20.8 & 23.4 & 13.3 \\ 
        X-n228-k23 & 10379.97 & 10379.84 & 21.8 & 239.4 & 10.2 \\ 
        X-n233-k16 & 8480.25 & 8479.83 & 19.6 & 253.3 & 6.3 \\ 
        X-n237-k14 & 10870.43 & 10870.43 & 21.7 & 34.9 & 5.0 \\ 
        X-n242-k48 & 28644.96 & 28640.61 & 27.3 & 123.6 & 17.2 \\ \hline
 \end{tabular}
 }
\end{table}

\begin{table}[htbp]
 \fontsize{10}{12}\selectfont
\caption{Performance of the route-optimization step, for MVRPB instances with 100 clients per period} \label{tab3}
 \centering
 \setlength\tabcolsep{0.15cm}
 \scalebox{0.95}
 {
 \begin{tabular}{c@{\hspace*{1cm}}rr@{\hspace*{1cm}}rr@{\hspace*{1cm}}r}
 \hline
               \multirow{2}{*}{Instance} & \multicolumn{2}{c}{Distance} & \multicolumn{2}{c}{T(s)} &\multirow{2}{*}{\#k}  \\ 
           \cmidrule(lr){2-3}\cmidrule(lr){4-5}
       & HGS & VRPSolver & HGS & VRPSolver & \\ \hline
        X-n200-k36 & 30780.5 & 30758.06 & 47.6 & 979.1 & 19.6 \\ 
        X-n204-k19 & 11535.58 & 11534.75 & 32.5 & 470.4 & 10.0 \\ 
        X-n209-k16 & 16724.49 & 16722.15 & 39.9 & 420.6 & 8.3 \\ 
        X-n214-k11 & 6076.82 & 6076.64 & 39.2 & 1041.6 & 6.1 \\ 
        X-n219-k73 & 55473.58 & 55473.58 & 27.6 & 2.7 & 34.0 \\ 
        X-n223-k34 & 20522.69 & 20519.9 & 31.8 & 69.3 & 17.1 \\ 
        X-n228-k23 & 12888.06 & 12887.6 & 32.2 & 950.6 & 12.5 \\ 
        X-n233-k16 & 10288.46 & 10287.92 &27.0 & 1296.1 & 8.0 \\ 
        X-n237-k14 & 13518.71 & 13518.58 & 35.0 & 697.4 & 6.0 \\ 
        X-n242-k48 & 37262.13 & 37247.82 & 40.9 & 324.4 & 22.4 \\ \hline
 \end{tabular}
 }
\end{table}

As seen in these experiments, the computational time used by VRPSolver to optimally solve the underlying CVRP problem for each period is generally small, with a median value of $25.0$ seconds. However, on a handful of cases, the computational time may be long, reaching $12.8$ hours in the worst case (one exceptional case in 3\,000 CVRP single-period sub-problems). In contrast, HGS has a more controllable computational time, ranging from $11.4$ seconds to $2.43$ minutes, with a median value of $21.8$ seconds. We observe that the initial solutions found by HGS were almost optimal in terms of their distance, with an average distance over all instances of $16715.0$ compared to $16712.9$ for VRPSolver, i.e., with an average gap error of only $0.013 \%$ from optimal solution values. Given this, we recommend using HGS as the underlying solution approach for the route optimization step in practical time-critical applications. In the context of this study, we decided to complete the solution process to achieve proven optima with VRPSolver, as this will subsequently permit us to derive bounds on the best possible workload balance through Equation~\eqref{eq:optbound}.

Finally, we must observe that the problems associated with each period are independent, such that it is possible to solve them in parallel. We used this observation in our experiments, as the multi-core structure of our processor permitted us to independently solve the CVRPs associated with each period on a different core, therefore maximizing our utilization of available computational resources and reducing the total time needed to conduct our experiments.

\noindent
\textbf{Multi-period workload balancing.}
In the workload balancing step, the routes of the optimal CVRP solution for each period are assigned to drivers to minimize the maximum (\textsc{Min-Max}) total distance traveled by each driver over the entire planning horizon. We build our analysis on three key workload measurements:
\begin{itemize}[nosep]
    \item UB -- The initial workload produced by the constructive approach described in Section~\ref{UB}. The workload corresponds to the largest total distance for a driver over the entire planning horizon.
    \item Opt -- The optimal workload obtained after completing the binary search.
    \item LB -- The lower bound of Equation~\eqref{eq:optbound}, which assumed that distance is optimal and workload equity is perfect (often this does not match a practical solution).
\end{itemize}

Tables~\ref{tab50-2:wmax-new}~to~\ref{tab100-2:wmax-new} report the workload values of UB, Opt, and LB for the different instances, with a varying number of periods and with 50, 75, and 100 clients in each period, respectively. 
Each line in the tables corresponds to an average value over ten different MVRPB instances. These tables also indicate the average number of binary-search operations in our algorithm and the average solution time.

\begin{landscape}
\thispagestyle{empty}
\begin{table}[H]
\centering
\vspace*{-1cm}
\caption{Performance of the multi-period workload balancing step -- MVRPB instances with 50 clients per period} 
\label{tab50-2:wmax-new}
\vspace*{-0.3cm}
\setlength\tabcolsep{0.15cm}
\hspace*{-1cm}
\scalebox{0.6}
 {\begin{tabular}{c@{\hspace*{1cm}}rrrrr@{\hspace*{1cm}}rrrrr@{\hspace*{1cm}}rrrrr@{\hspace*{1cm}}rrrrr@{\hspace*{1cm}}rrrrr}
    \toprule
        \multirow{2}{*}{Tests} & \multicolumn{5}{c}{2 periods}  & \multicolumn{5}{c}{3 periods} & \multicolumn{5}{c}{5 periods} & \multicolumn{5}{c}{7 periods} &   \multicolumn{5}{c}{10 periods}  \\ 
        \cmidrule(lr){2-6}\cmidrule(lr){7-11}\cmidrule(lr){12-16}\cmidrule(lr){17-21}\cmidrule(lr){22-26}
          &LB & UB & Opt & \#It & T(s) & LB & UB & Opt & \#It & T(s) & LB & UB & Opt & \#It & T(s) & LB & UB & Opt & \#It & T(s) & LB & UB & Opt & \#It & T(s) \\ \hline
        X-n200-k36 & 3199.6 & 3551.5 & 3551.5 & 10.0 & 6.46 & 4858.2 & 5291.1 & 5033.9 & 9.7 & 7.18 & 8084.7 & 8630.9 & 8090.3 & 10.0 & 56.35 & 11372.8 & 11989.7 & 11372.8 & 1.0 & 8.47 & 16274.3 & 16907.1 & 16274.5 & 3.0 & 279.05 \\ 
        X-n204-k19 & 2735.3 & 3185.4 & 3185.4 & 10.1 & 6.54 & 4107.5 & 4570.4 & 4320.7 & 9.8 & 6.80 & 6810.6 & 7095.7 & 6822.5 & 9.1 & 5.82 & 9584.2 & 9721.3 & 9585.0 & 7.0 & 18.75 & 13663.1 & 13883.5 & 13663.1 & 1.0 & 9.39 \\
        X-n209-k16 & 4077.9 & 4846.0 & 4846.0 & 11.0 & 6.65 & 6167.3 & 7088.4 & 6421.1 & 10.9 & 7.36 & 10276.2 & 10574.6 & 10304.9 & 9.1 & 5.78 & 14297.5 & 14690.3 & 14300.5 & 9.4 & 23.59 & 20353.8 & 20725.2 & 20353.8 & 1.0 & 8.15 \\ 
        X-n214-k11 & 2256.0 & 2644.0 & 2644.0 & 10.1 & 6.39 & 3403.7 & 3643.4 & 3546.0 & 8.9 & 5.98 & 5608.7 & 5886.0 & 5652.8 & 9.2 & 5.89 & 7939.2 & 8186.8 & 7946.5 & 8.7 & 14.93 & 11333.4 & 11498.9 & 11333.5 & 1.6 & 8.62 \\ 
        X-n219-k73 & 3378.8 & 3696.6 & 3696.6 & 9.9 & 6.68 & 5012.3 & 5469.3 & 5035.8 & 9.9 & 6.03 & 8377.7 & 8650.5 & 8377.7 & 1.0 & 740.10 & 11778.3 & 12115.5 & 11778.3 & 2.2 & 17.72 & 16801.8 & 17059.7 & 16801.8 & 1.0 & 23.96 \\ 
        X-n223-k34 & 2436.9 & 2909.9 & 2909.9 & 10.3 & 6.19 & 3630.5 & 4052.6 & 3708.1 & 9.6 & 5.34 & 5991.6 & 6285.4 & 5994.4 & 9.3 & 49.34 & 8438.0 & 8702.3 & 8438.0 & 1.0 & 7.77 & 12101.1 & 12276.0 & 12101.2 & 1.7 & 439.17 \\ 
        X-n228-k23 & 2084.1 & 2705.4 & 2705.4 & 10.7 & 6.30 & 3156.1 & 3665.2 & 3412.0 & 10.0 & 4.98 & 5254.3 & 5650.0 & 5271.0 & 9.3 & 12.61 & 7291.4 & 7592.9 & 7291.5 & 1.9 & 11.91 & 10513.6 & 10815.5 & 10513.6 & 1.3 & 9.87 \\ 
        X-n233-k16 & 2998.2 & 3541.6 & 3541.6 & 10.5 & 6.73 & 4524.8 & 5116.2 & 4785.1 & 9.9 & 6.18 & 7560.4 & 8115.5 & 7594.1 & 9.7 & 6.62 & 10527.6 & 11028.3 & 10530.8 & 9.8 & 20.70 & 15072.1 & 15367.5 & 15072.1 & 1.0 & 7.69 \\ 
       X-n237-k14 & 5336.4 & 5718.7 & 5718.7 & 9.7 & 6.39 & 7995.7 & 8451.0 & 8196.9 & 9.9 & 5.63 & 13326.9 & 13647.9 & 13384.9 & 8.4 & 6.71 & 18636.1 & 18856.6 & 18646.0 & 8.2 & 8.15 & 26601.9 & 26985.1 & 26602.1 & 2.6 & 8.58 \\ 
        X-n242-k48 & 3329.5 & 3875.0 & 3875.0 & 10.6 & 6.71 & 4957.2 & 5583.6 & 5063.6 & 10.2 & 5.53 & 8224.7 & 8669.3 & 8226.3 & 9.7 & 52.25 & 11365.1 & 11899.5 & 11509.3 & 2.1 & 10.54 & 16527.5 & 16799.3 & 16527.5 & 1.0 & 18.65 \\ \hline
        Ave. & 3183.3 & 3667.4 & 3667.4 & 10.3 & 6.50 & 4781.3 & 5293.1 & 4952.3 & 9.9 & 6.10 & 7951.6 & 8320.6 & 7971.9 & 8.5 & 94.15 & 11123.0 & 11478.3 & 11139.9 & 5.1 & 14.25 & 15924.3 & 16231.8 & 15924.3 & 1.5 & 81.31 \\ \hline
    \end{tabular}}
\end{table}

\begin{table}[H]
\centering
\vspace*{-0.45cm}
\caption{Performance of the multi-period workload balancing step -- MVRPB instances with 75 clients per period} 
\label{tab75-2:wmax-new}
\vspace*{-0.3cm}
\setlength\tabcolsep{0.15cm}
\hspace*{-1cm}
\scalebox{0.6}
 {\begin{tabular}{c@{\hspace*{1cm}}rrrrr@{\hspace*{1cm}}rrrrr@{\hspace*{1cm}}rrrrr@{\hspace*{1cm}}rrrrr@{\hspace*{1cm}}rrrrr}
    \toprule
        \multirow{2}{*}{Tests} & \multicolumn{5}{c}{2 periods}  & \multicolumn{5}{c}{3 periods} & \multicolumn{5}{c}{5 periods} & \multicolumn{5}{c}{7 periods} &   \multicolumn{5}{c}{10 periods}  \\ 
        \cmidrule(lr){2-6}\cmidrule(lr){7-11}\cmidrule(lr){12-16}\cmidrule(lr){17-21}\cmidrule(lr){22-26}
         & LB & UB & Opt & \#It & T(s) & LB & UB & Opt & \#It & T(s) & LB & UB & Opt & \#It & T(s) & LB & UB & Opt & \#It & T(s) & LB & UB & Opt & \#It & T(s) \\ \hline
        X-n200-k36 & 3201.3 & 3463.1 & 3463.1 & 9.6 & 6.67 & 4788.9 & 5201.0 & 4916.3 & 9.8 & 7.60 & 7954.8 & 8459.4 & 7955.3 & 6.0 & 23.64 & 11132.6 & 11608.2 & 11132.6 & 1.0 & 11.97 & 15919.8 & 16472.8 & 15919.8 & 1.0 & 20.89 \\ 
        X-n204-k19 & 2345.8 & 2695.6 & 2695.6 & 9.8 & 6.78 & 3540.7 & 3841.5 & 3622.5 & 9.4 & 6.78 & 5895.7 & 6219.6 & 5898.7 & 9.6 & 27.47 & 8257.0 & 8633.6 & 8257.0 & 1.0 & 7.93 & 11786.4 & 11941.7 & 11786.4 & 1.0 & 13.86 \\
        X-n209-k16 & 4089.0 & 4606.9 & 4606.9 & 10.4 & 6.62 & 6136.5 & 6724.1 & 6284.7 & 10.1 & 6.68 & 10188.3 & 10686.7 & 10196.7 & 9.9 & 16.26 & 14216.4 & 14493.9 & 14216.8 & 4.2 & 22.23 & 20224.5 & 20496.7 & 20224.5 & 1.0 & 15.97 \\
        X-n214-k11 & 1951.2 & 2292.2 & 2292.2 & 9.7 & 6.29 & 2918.5 & 3289.1 & 3090.3 & 9.7 & 6.40 & 4846.5 & 5115.8 & 4859.4 & 9.1 & 8.21 & 6805.0 & 7026.2 & 6805.5 & 5.1 & 13.41 & 9762.8 & 9945.0 & 9762.8 & 1.0 & 9.21 \\ 
        X-n219-k73 & 3349.5 & 3604.7 & 3604.7 & 9.6 & 7.28 & 5002.9 & 5498.0 & 5014.2 & 9.7 & 11.77 & 8337.8 & 8553.9 & 8337.9 & 1.8 & 1705.40 & 11696.1 & 11979.5 & 11696.1 & 1.0 & 20.36 & 16731.2 & 17034.0 & 16731.3 & 1.9 & 52.37 \\ 
        X-n223-k34 & 2424.1 & 2770.6 & 2770.6 & 9.9 & 6.71 & 3604.4 & 3978.5 & 3653.8 & 9.5 & 7.57 & 5983.5 & 6219.6 & 5984.0 & 5.2 & 132.47 & 8404.2 & 8601.1 & 8404.2 & 1.0 & 10.84 & 12071.3 & 12249.7 & 12071.3 & 1.0 & 22.66 \\ 
        X-n228-k23 & 2033.6 & 2466.4 & 2462.6 & 10.0 & 6.82 & 3053.9 & 3622.9 & 3188.7 & 10.1 & 7.28 & 5017.0 & 5348.4 & 5020.6 & 9.7 & 51.00 & 7103.5 & 7431.2 & 7103.5 & 1.0 & 7.72 & 10217.8 & 10484.5 & 10217.8 & 1.0 & 17.36 \\ 
        X-n233-k16 & 2709.8 & 3103.9 & 3103.9 & 10.0 & 6.62 & 4055.7 & 4570.3 & 4251.0 & 9.9 & 7.38 & 6818.4 & 7142.7 & 6827.9 & 9.5 & 11.17 & 9532.3 & 9847.3 & 9532.4 & 1.9 & 11.33 & 13525.2 & 13792.1 & 13525.2 & 1.0 & 13.70 \\ 
        X-n237-k14 & 4313.7 & 5148.3 & 5148.3 & 11.1 & 7.14 & 6472.5 & 7475.0 & 6965.4 & 10.9 & 7.04 & 10835.0 & 11452.6 & 10854.5 & 9.2 & 7.68 & 15187.6 & 15520.8 & 15189.3 & 9.4 & 20.20 & 21741.2 & 21826.9 & 21741.2 & 1.0 & 8.68 \\ 
        X-n242-k48 & 3411.4 & 3886.9 & 3886.9 & 10.4 & 6.97 & 5107.4 & 5707.5 & 5171.3 & 10.1 & 8.16 & 8470.0 & 8984.9 & 8470.2 & 2.8 & 194.67 & 11752.5 & 12041.8 & 11752.5 & 1.0 & 12.82 & 16658.9 & 16882.6 & 16658.9 & 1.0 & 22.91 \\ \hline
        Ave. & 2982.9 & 3403.9 & 3403.5 & 10.1 & 6.79 & 4468.1 & 4990.8 & 4615.8 & 9.9 & 7.67 & 7434.7 & 7818.4 & 7440.5 & 7.3 & 217.80 & 10408.7 & 10718.4 & 10409.0 & 2.7 & 13.88 & 14863.9 & 15112.6 & 14863.9 & 1.1 & 19.76 \\ \hline
    \end{tabular}}
\end{table}

\begin{table}[H]
\centering
\vspace*{-0.45cm}
\caption{Performance of the multi-period workload balancing step -- MVRPB instances with 100 clients per period} 
\label{tab100-2:wmax-new}
\vspace*{-0.3cm}
\setlength\tabcolsep{0.15cm}
\hspace*{-1cm}
\scalebox{0.6}
 {\begin{tabular}{c@{\hspace*{1cm}}rrrrr@{\hspace*{1cm}}rrrrr@{\hspace*{1cm}}rrrrr@{\hspace*{1cm}}rrrrr@{\hspace*{1cm}}rrrrr}
    \toprule
        \multirow{2}{*}{Tests} & \multicolumn{5}{c}{2 periods}  & \multicolumn{5}{c}{3 periods} & \multicolumn{5}{c}{5 periods} & \multicolumn{5}{c}{7 periods} &   \multicolumn{5}{c}{10 periods}  \\ 
        \cmidrule(lr){2-6}\cmidrule(lr){7-11}\cmidrule(lr){12-16}\cmidrule(lr){17-21}\cmidrule(lr){22-26}
         & LB & UB & Opt & \#It & T(s) & LB & UB & Opt & \#It & T(s) & LB & UB & Opt & \#It & T(s) & LB & UB & Opt & \#It & T(s) & LB & UB & Opt & \#It & T(s) \\ \hline
        X-n200-k36 & 3161.3 & 3439.1 & 3407.4 & 9.4 & 7.61 & 4727.0 & 5311.2 & 4834.0 & 10.0 & 8.86 & 7865.5 & 8477.5 & 7865.5 & 1.0 & 76.79 & 10982.5 & 11470.7 & 10982.5 & 1.0 & 13.99 & 15703.2 & 16244.7 & 15703.2 & 1.0 & 25.56 \\ 
        X-n204-k19 & 2309.4 & 2607.1 & 2607.1 & 9.6 & 7.21 & 3460.6 & 3731.8 & 3507.1 & 9.1 & 7.99 & 5763.0 & 6048.3 & 5764.1 & 9.5 & 90.53 & 8071.2 & 8337.0 & 8071.2 & 1.0 & 10.42 & 11527.4 & 11775.4 & 11527.4 & 1.0 & 12.66 \\ 
        X-n209-k16 & 4043.6 & 4448.2 & 4448.2 & 10.1 & 8.00 & 6062.4 & 6610.1 & 6173.2 & 10.1 & 6.33 & 10137.2 & 10570.4 & 10140.8 & 9.8 & 35.44 & 14176.1 & 14523.5 & 14176.1 & 1.0 & 8.28 & 20204.7 & 20512.9 & 20204.7 & 1.0 & 11.06 \\ 
        X-n214-k11 & 1969.8 & 2176.2 & 2176.2 & 9.2 & 7.96 & 2967.3 & 3311.4 & 3053.0 & 9.2 & 5.91 & 4987.2 & 5239.6 & 4994.1 & 8.9 & 10.83 & 7026.2 & 7299.6 & 7026.4 & 2.5 & 10.80 & 9986.3 & 10105.7 & 9986.3 & 1.0 & 9.27 \\
        X-n219-k73 & 3239.3 & 3486.3 & 3486.3 & 9.5 & 7.36 & 4879.9 & 5296.2 & 4886.1 & 9.8 & 21.87 & 8159.9 & 8477.6 & 8159.9 & 1.0 & 15.04 & 11397.0 & 11765.1 & 11397.0 & 1.0 & 31.32 & 16316.2 & 16654.4 & 16316.2 & 1.0 & 102.17 \\ 
        X-n223-k34 & 2384.1 & 2708.7 & 2708.7 & 9.7 & 6.95 & 3603.1 & 3871.1 & 3631.7 & 9.1 & 9.55 & 6054.6 & 6243.3 & 6054.6 & 1.0 & 433.94 & 8413.6 & 8631.7 & 8413.6 & 1.0 & 12.68 & 12001.9 & 12232.6 & 12001.9 & 1.0 & 20.40 \\ 
        X-n228-k23 & 2012.2 & 2479.0 & 2479.0 & 10.3 & 7.32 & 3064.6 & 3598.2 & 3180.6 & 10.0 & 8.19 & 5112.6 & 5417.1 & 5114.5 & 7.4 & 29.32 & 7169.0 & 7410.0 & 7169.0 & 1.0 & 9.67 & 10325.3 & 10654.1 & 10325.3 & 1.0 & 12.57 \\
        X-n233-k16 & 2524.4 & 2955.2 & 2923.2 & 10.0 & 7.58 & 3754.3 & 4308.4 & 3885.5 & 10.0 & 7.79 & 6307.1 & 6874.5 & 6310.2 & 9.9 & 28.01 & 8849.7 & 9410.8 & 8849.7 & 1.0 & 11.70 & 12602.4 & 12996.6 & 12602.4 & 1.0 & 10.34 \\ 
        X-n237-k14 & 4487.4 & 4942.2 & 4942.2 & 10.3 & 7.41 & 6733.5 & 7343.4 & 6919.2 & 10.3 & 7.22 & 11209.6 & 11801.3 & 11218.6 & 10.1 & 11.40 & 15710.0 & 15927.3 & 15710.3 & 3.4 & 14.09 & 22531.5 & 22764.5 & 22531.5 & 1.0 & 9.66 \\ 
        X-n242-k48 & 3386.4 & 3832.8 & 3832.8 & 10.3 & 7.17 & 5049.7 & 5545.3 & 5072.4 & 9.8 & 10.15 & 8254.0 & 8804.4 & 8254.0 & 1.0 & 1467.53 & 11637.7 & 12000.6 & 11637.9 & 2.9 & 552.76 & 16633.5 & 16834.3 & 16633.5 & 1.0 & 32.24 \\ \hline
        Ave. & 2951.8 & 3307.5 & 3301.1 & 9.8 & 7.46 & 4430.2 & 4892.7 & 4514.3 & 9.7 & 9.39 & 7385.1 & 7795.4 & 7387.6 & 6.0 & 219.88 & 10343.3 & 10677.6 & 10343.4 & 1.6 & 67.57 & 14783.2 & 15077.5 & 14783.2 & 1.0 & 24.59 \\ \hline
    \end{tabular}}
\end{table}
\end{landscape}

As seen in these experiments, only a few seconds are required in most cases to complete an optimal multi-period workload balancing step. This confirms the efficiency of our two-stage solution approach. Generally, instances involving a larger number of routes and drivers (e.g., X-n219-k73) lead to more complex workload balancing problems. Again, there is also inherent variability due to the exact solution process, given that MILP approaches can exhibit substantially different computation times when solving instances of similar sizes. Overall, the computational time of the second step ranges from 4.40 to 7600.99 seconds with a median value of 8.53 seconds.

The workload allocation created by the initial construction approach (described in Section~\ref{UB}) is optimal (i.e., UB = Opt) for 298 out of 300 MVRPB instances with 2 periods, as well as for 15 out of 300 MVRPB instances with 3 periods. In contrast, as soon as the number of periods becomes greater than five, the initial construction approach is unlikely to lead to the best possible workload allocation, and the underlying mathematical model produces much better solutions.

When the planning horizon contains five days or more, we frequently notice cases of ``perfect'' workload equity, where the obtained workload balance matches the theoretical lower bound (i.e., Opt = LB). The ability to achieve perfect balance comes as a consequence of the increased number of possible assignment solutions, which grows exponentially with the number of periods. In practice, in cases with five periods or more, it is sufficient to focus on cost-optimal routing solutions and to create equitable workloads by careful assignment. Finally, as our binary search strategy includes a first step to verify if a solution with perfect balance exists, all the cases for which perfect balance is possible are solved in a single call to the feasibility subproblem. In the other cases, it generally takes between $7$ to $11$ steps.

\subsection{Planning-Horizon Length and Workload Equity} 
\label{EMC}

\begin{figure}[htbp]
\begin{center}
\centering
\includegraphics[width=0.83\textwidth]{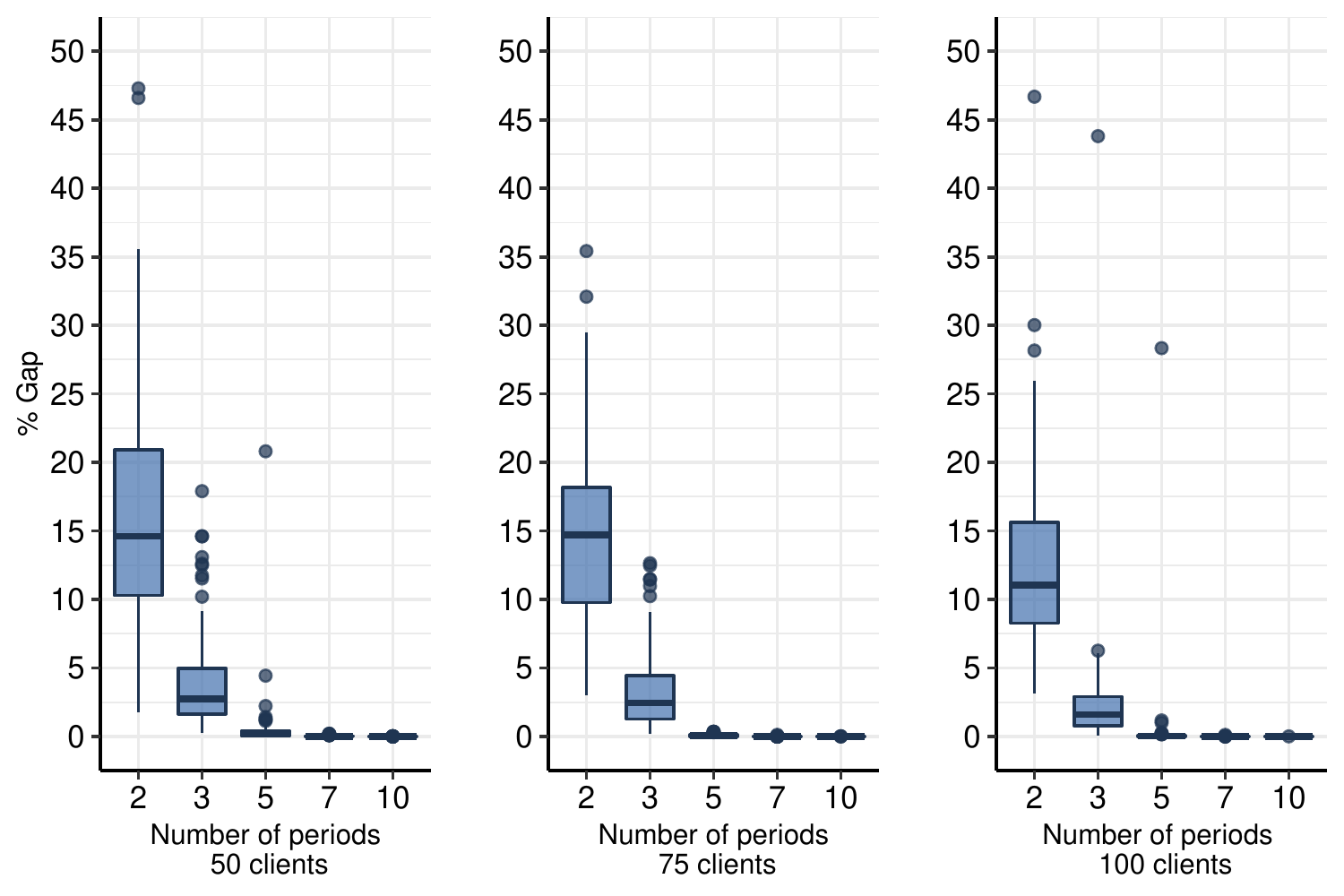}
\caption{Convergence of Opt toward LB as the number of periods increases for instances with 50, 75, and 100 clients in each period.}
\label{deviation}
\end{center}
\end{figure}

The previous section showed that near-perfect workload equity is achievable, in practice, for planning horizons with at least five periods. To visualize more clearly the impact of the number of periods over workload equity, Figure~\ref{deviation} provides additional boxplot representations of the Gap(\%) between the ideal workload (LB from Equation~\eqref{eq:optbound}) and the optimal solution value (Opt) of our two-phases approach. We used the following calculation: $\text{Gap} = 100 \times (\text{Opt}-\text{LB})/\text{LB}$. We provide separate plots for the cases with 50, 75, or 100 customers per period. Each boxplot corresponds to the data of a given planning horizon with $T \in \{2, 3, 5, 7, 10\}$ periods, therefore gathering gap measurements from 100 instances. The whiskers indicate the minimum, first quartile (Q1), median, third quartile (Q3), and maximum. The minimum corresponds to $Q1 - 1.5 \times \text{interquartile range}$, while the maximum corresponds to $Q3 + 1.5 \times \text{interquartile range}$. Outliers that fall beyond the minimum and maximum range are additionally depicted as small circles.

These boxplots give another viewpoint on the convergence toward the best possible equity as the number of periods increases. In the vast majority of the cases (excluding a few outliers), a planning horizon of five days is sufficient to find equitable solutions with workload discrepancies below 1\% between drivers. In these situations, there is no need to seek a trade-off between routing costs and workload equity since optimal routing solutions can be used to achieve equity. Another benefit of our two-stage approach is its flexibility since additional constraints, decisions, and objectives (i.e.,  routing attributes -- \citealt{VIDAL20131}) only need to be integrated with the first phase of the solution approach.

Finally, in the cases with very few periods (e.g., 2 or 3 days), we observe that focusing the search on optimal routing solutions does not permit achieving the best possible workload equity. In such situations, it would be helpful to consider alternative routing solutions. One possibility would consist in producing multiple routing solutions for each period and extending the workload balancing step to include all these alternatives. Another approach, more complex to develop in practice, would be to solve the routing and driver-allocation problem in an integrated manner, considering distance and workload equity in a bi-objective solution method. However, in both cases, the user would need to specify a trade-off between acceptable extra routing costs and the desired workload equity level.

\section{Conclusions}
\label{conclusion}

In this work, we have revisited workload equity in vehicle routing with a longer-term perspective, considering a planning horizon of several days. We have shown that a two-phase optimization approach can identify the most equitable solutions with minimal distance. When the planning horizon exceeds five days, the resulting solutions are optimal in terms of distance and near-optimal (below 1\% gap) in terms of equity. Therefore, workload equity appears to be achievable without integrated approaches and trade-off calibration, and without any compromise on operational efficiency.

Several important research perspectives are open in connection with this study. Firstly, our work focused on deterministic settings, where the complete customer demand is known on the planning horizon. Practical situations often involve dynamically-revealed request information, and therefore it is an open question to determine to what extent multi-period workload equity is achievable in dynamic contexts. Another important aspect of practical delivery systems concerns delivery consistency. When the same driver regularly visits the same areas or clients, the service quality and the satisfaction of drivers and clients generally increases \citep{articleKovacs}. Our approach towards equitable workload allocation largely benefits from the exponential number of possible route-driver allocation combinations. However, consistency may significantly reduce the number of allocation possibilities, such that new approaches may be needed to conciliate three key aspects in a multi-period setting: cost efficiency, workload equity, and delivery consistency.\\

\noindent\textbf{Acknowledgments.}
Financial support for this research has been provided by the Natural Sciences and Engineering Research Council of Canada (NSERC) and by an IVADO Fundamental Research grant. This study has been also partly enabled by the computational infrastructure provided by Calcul Québec and Compute Canada. This support is gratefully acknowledged.

\end{document}